\documentclass[12pt]{article}
\usepackage[margin=0.8in]{geometry}
\usepackage{graphicx,amsthm,amsmath,amsfonts,latexsym,amssymb,multirow,booktabs}
\usepackage{xcolor,array,verbatim,bm}
\usepackage{xspace}

\usepackage[T1]{fontenc}
\usepackage{mathpazo} 
\usepackage{tgpagella} 

\definecolor{linkbrown}{rgb}{0.50,0.22,0.10}
\usepackage[round]{natbib}
\usepackage[colorlinks=true,linkcolor=linkbrown,citecolor=linkbrown,urlcolor=linkbrown]{hyperref}
\title{The Pirah\~a and the cognitive gap in Frege's theorem: Hume's principle without the \#}
\author{Subrata~Pal\\
{\small Department of Neurology, Washington University in St Louis.}}

\newcommand{\Piraha}{Pirah\~a\xspace}
\newtheorem{principle}{Principle}

\begin{document}
 \maketitle
\begin{abstract}
Frege's theorem proves that Hume's principle, in second-order logic, yields all of arithmetic. Yet the \Piraha{} people show one-to-one correspondence (equinumerosity) only where pairing can be enacted, with its range extended under local training, and still have no counting or arithmetic. We argue this is not a paradox but a matter of precise localization. Hume's principle includes a \emph{cardinality operator \#} that names cardinals as objects (often modeled as equivalence classes of equinumerous concepts), and what the \Piraha{} lack is not the relation but this operator. We identify the \emph{number-word practice} as the \emph{cognitive realization} of \#, which recasts the ``number-as-cognitive-technology'' thesis in formal terms and locates the cognitive boundary at symbolization, not recursion.

The identification is generative, not decorative: the reach of \# tracks the reach of the token practice that carries it, so across languages and cultures we see a gradient, not a sharp cliff. 
And number words are not special as words; what \# needs is any stable, reusable marker that can preserve exact cardinal identity across absence, rearrangement, delay, or modality shift: a spoken numeral, a scratch on a stick, or a knot in a cord. So the thesis is about having some symbolic token-practice, not about language specifically.
It is supported by converging evidence from Nicaraguan homesigners, numerate adults under verbal interference, and cross-linguistic numeral gradients. We make no causal, acquisition, or neural claim; the identification is constitutive.
\end{abstract}

\noindent%

\section{Introduction}

Daniel L. Everett~\footnote{Keren Madora (n\'ee Graham), Daniel Everett's former wife, is the only non-\Piraha to have lived longer among them than he did \citep{everett2005}. She is coauthor of a later replication \citep{everettandmadora2012}; Section~\ref{sec:about-Pirahas} reconciles its apparent disagreement with \citet{franketal2008}.} went to the Maici river in 1977 \citep{everett1986} as a missionary and stayed to become one of the principal anthropological linguists of the \Piraha{}\footnote{Throughout, `\Piraha' names both the people (a few hundred hunter-gatherers on the Maici, a tributary of the Madeira; population estimates range from roughly 200 to 700 across sources) and their language, the last surviving member of the Mura family. The people's own name for themselves is hiait\'{\i}ih\'{\i} (`a straight one'). They call their language 'apait\'{\i}iso (`a straight head') and all other languages 'apag\'aiso (`a crooked head') \citep{everett2005}.}. The language he described is rich, with verbs carrying up to sixteen optional suffix slots, among them evidential markers for how the speaker knows what he asserts \citep{everett2005}. However, it became famous for what it lacks: Everett reported that \Piraha has no embedded clauses. It ignited one of the loudest disputes in modern linguistics as it directly attacks the concept of subordinate clauses and recursion being the defining core of the human language faculty \citep{hauseretal2002, everett2005, nevinsetal2009, everett2009}. This paper takes no side in it.

Our subject is a quieter absence in the same language which no party to the dispute denies: \Piraha has no numerals, not even a stable word for one \citep{gordon2004, everett2005, franketal2008}. The absence mattered to the \Piraha{} themselves. 
The \Piraha{} asked Daniel and Keren to teach them number systems \citep{everett2005}, so that other people would not cheat them. Eight months of evening classes produced no one who could count to ten or add one and one.

Why should it be so hard? We made a natural first conjecture: \emph{counting seems to need a successor step, and successor-like structure is what subordination supplies, so perhaps a language without embedding cannot host counting}. We will argue the break lies elsewhere.

Now we step back and look at counting through the philosopher's glass: what, minimally, must a mind possess for arithmetic?
For the axiomatization of arithmetic, Leibniz started by defining `one' and `increase by one' \citep[IV.vii.10]{leibniz1704}. Later, Peano, Dedekind, Cantor, Frege, and their successors tried to develop such axioms in a consistent way. Specifically, Peano and Dedekind's approach \citep{peano1889} is prevalent in the mainstream mathematics community (see \citealp{tao2006} for a clean modern statement), where a successor function is essential for creating the whole number system.
Frege (with later reconstructions by Wright, Boolos, and Heck; Section~\ref{sec:HP}) came from a slightly different direction. From that direction it can be shown that second-order logic and a single statement, Hume's principle (HP), together yield all of arithmetic (precisely, the Dedekind-Peano axioms follow as theorems). Hume's principle is based on a seemingly inconsequential part of the book \emph{A Treatise of Human Nature}, which says: ``When two \emph{numbers} [sets, in the usage of the time] are so combin'd, as that the one has always an unit 
answering to every unit of the other, we pronounce them equal'' \citep[I.iii.1]{hume1739}.
In the modern paraphrase it is deceptively simple. Let \emph{cardinality} mean the sheer how-many of a collection, and let F, G be any two collections: ``the cardinality of F and the cardinality of G are the same if, and only if, F and G are equinumerous.''

Return now to the \Piraha. 
In \citeauthor{gordon2004}'s matching tasks \citeyearpar{gordon2004}, \Piraha{} adults placed objects one against each target, though with notable errors at larger set sizes. \citet{franketal2008} repeated the task with a more controlled design and found the same one-to-one matching essentially without error, out to ten objects. The two results disagree; Section~\ref{sec:about-Pirahas} names the hidden variable that reconciles them.
The picture is not unique to the Amazon: in Nicaragua, deaf homesigners, adults who grew up inventing their own gestures inside a fully numerate culture, have no count list either; one of them, describing a set of ten sheep, extends nine fingers \citep{spaepenetal2011}. We will return to this image (Section~\ref{sec:scopes}).

The two stories now meet in a single question: ``If Hume's principle is all arithmetic needs, and the \Piraha{} hold its equinumerosity side in the enacted, finite sense below, why does matching not deliver counting?''

We investigate this question formally and show that the gap is precisely the cardinality operator \#, and \textbf{the number-word practice is the cognitive realization of that operator}, which \citet{franketal2008} call \emph{cognitive technology}. 

It is important to notice that, i) we do not have any new data and reinterpret published findings only; ii) we do not form any strong-Whorfian claim and in fact argue that such claims are statistically ill-formed given the available data (the candidate variables are confounded; Section~\ref{sec:scopes}); iii) we do not claim anything about whether \Piraha has syntactic recursion; our argument is independent of that question; iv) we make no claim about acquisition either.

\section{Two sides in detail}
\subsection{The empirical side}\label{sec:about-Pirahas}
It was well known that the \Piraha{} do not have a discrete number system. The language has three quantity words, \emph{h\'oi}, \emph{ho\'i}, and \emph{ba\'agiso} (roughly `a small amount,' `a larger amount,' and `many', respectively).  
\citet{gordon2004} devised a test which tried to capture whether they have one-to-one matching ability beyond these three \emph{concept-words}. 
However, their study had only seven participants (six men and one woman; most data from four of the men)\footnote{Males trade and have more contact with Portuguese and other languages, so one could imagine an individual exposure gradient. It does not fit: \citet{franketal2008}'s sex-balanced sample sat at ceiling on the two direct-pairing tasks, and the equally sex-balanced sample of \citet{everettandmadora2012} did not. The operative exposure is village-level (see below), not individual. No by-sex breakdown is published, so this stays a footnote, not a battle.} and was followed by \citet{franketal2008} with a larger, sex-balanced sample (n=14) and a more controlled design. The first design elicited what those three words mean, in both increasing and decreasing order of quantity, which suggested that those words are approximate. The second design, independent of the first, consisted of five subtasks of increasing difficulty, each asking the participant to put out exactly as many balloons as there were spools. In the first two (one-to-one line match, uneven match), each balloon can be laid directly against its spool, and the error is essentially nil. In the remaining three (orthogonal match, hidden match, nuts-in-a-can), no direct placement is possible: the quantity itself has to be carried across a rotation, a screen, or time, and there the errors are much higher. Note that in the orthogonal match the spools stay in full view; what is removed is only the possibility of pairing directly against them.

\citet{everettandmadora2012} re-tested the same one-to-one matching and found the \Piraha{} fail it beyond three, only 2 of 14 error-free where \citet{franketal2008} had 13 of 14, which reads at first as a flat contradiction. It resolves once \emph{a hidden variable is named}. Frank tested Xagiopai, the single village where \Piraha{} number words for four to ten had been introduced by Keren Madora and the matching task drilled, amid a decade of unusually heavy outside contact. The villages without that exposure, \citeauthor{gordon2004}'s among them, stop at three and estimate above it. The two studies are then two arms of a single \emph{local-training contrast}. And the response is task-selective: the enacted match goes from 32 of 56 correct trials untrained to 54 of 56 trained, while the two carried control tasks (hidden and orthogonal) sit at 24 of 56 in both studies. Whatever the active ingredient (labels, drilling, outside contact, or their combination; Section~\ref{sec:scopes}), the training stretched enacted matching and moved nothing that has to be carried. So their very disagreement already shows the pattern our thesis will explain: exact matching reaches as far as the practice that supports it, and no further. Section~\ref{sec:scopes} takes this up as a within-\Piraha{} de-confounder.

However, all three studies consistently show, on the failing tasks, a \emph{near-constant ($\sim$0.15) coefficient of variation} (CV) of errors. And here comes the surprise: the CV matches quite well with similar tests on numerate adults when the experimental design blocks exact token counting \citep{whalenetal1999, cordesetal2001}. Counting aloud produces rare off-by-one errors whose CV decreases with set size (rate: $1/\sqrt{n}$, binomial with fixed slip probability); counting suppressed gives constant CV near 0.15. This similar ratio-scaled error suggests that the \Piraha{} have a clear approximate number system (ANS), an innate, approximate sense of quantity, without a discrete label. In plain terms, their errors look just like a numerate person's when counting aloud is blocked: the analog magnitude sense is intact, only the exact tokens are missing. There have been cross-cultural CV comparisons (Munduruk\'u vs French) \citep{picaetal2004}, a developmental trajectory study \citep{halberdaandfeigenson2008}, and a study showing that education sharpens the CV only modestly \citep{piazzaetal2013}.

However, no amount of CV-sharpening yields exactness: an estimate whose error shrinks is still an estimate. Exactness is a change of representational type, a discrete token in place of an analog magnitude, not the limit of an ever-sharper approximation.

A word on `capacity,' since the argument leans on it. We do not mean the empty sense in which anyone could be trained to anything. We mean the measured baseline just established: exact correspondence for the few objects the un-augmented mind takes in at a glance, analog estimation (that same $\sim$0.15 CV) for everything above. An image worth keeping: three apples can be carried in bare hands; fifty need a bag. Watch someone without a bag drop the fourth apple: you have learned where their hands stop, not that they lack hands. `Capacity' names the hands. `Reach' names how far a practice can extend those hands: a drilled label may extend reach, as at the trained village above. A full count routine does not add hands; it supplies the bag. And what the bag adds is not merely more reach but \emph{portability}: the how-many can leave the apples behind.

Finally, the dissociation is not tied to this one culture: verbal interference in numerate adults \citep{franketal2012}, the Nicaraguan homesigners \citep{spaepenetal2011}, and cross-linguistic lexicon variation \citep{picaetal2004} reproduce it under three independent designs; Section~\ref{sec:scopes} does the careful accounting.

\subsection{The formal side: Hume's principle (HP)/Frege's theorem}\label{sec:HP}
Frege's original work \citep{frege1884, frege1893} proposed that his Basic Law V (BLV) and second-order logic (SOL) generate the whole of arithmetic; \citet{russell1903} proved the system inconsistent (Russell's paradox). 
Much later, \citet{parsons1965} noticed that a different abstraction principle of the same form as BLV, taken with SOL, suffices to generate the Peano axioms (PA) without any inconsistency \citep[see also Burgess, Hazen, and Hodes, independently]{boolos1987}. \citet{wright1983} rediscovered the point and made it a program; \citet{boolos1990} proved it rigorously. (The ``without any inconsistency'' clause does real work: an inconsistent principle, BLV included, is ``sufficient'' for everything.)

\citet{heck1993} went through the whole of Frege's articles and showed precisely that in those proofs, BLV enters only to obtain HP and the arithmetic follows from there. A very good introductory description of the whole point is provided in \citet{heck1999, heck2011}, and following him, we will call the result that HP plus SOL entails PA \emph{Frege's Theorem}. 

Precisely, following the notation of \citet{frege1884, wright1983, zalta2023}, we define a concept, which is a predicate: a function from objects to truth-values. `\underline{\ \ \ } is a card on this table' is a concept. It is true of each card on the table, say the queen of hearts if she is there, and false for everything not satisfying that predicate. In set-theoretic terms, a concept corresponds to a set (its extension or characteristic function), but Frege works with concepts directly. 

Notice that the operator \# takes a concept F and returns an object, \#F, the number of things falling under F. Then we formally write the HP:
\begin{principle}[Hume's Principle]\label{HP}
    For two concepts F and G, \# F = \# G, iff, F is equinumerous with G.
\end{principle}

Note that HP has two sides: the right-hand side is a relation (F and G can be put in one-to-one correspondence), while the left-hand side is an identity between objects (\#F = \#G). The operator \# is what bridges them. Forming \#F names the cardinal as an object (in the class model, the equivalence class of all concepts equinumerous with F, the quotient of the relation); the \Piraha{}, we will argue, hold the relation but never form this quotient.
The two sides are easy to conflate, even in print: \citet{izardetal2014}, studying one-to-one correspondence in 3-year-olds, give as ``Hume's principle'' the correspondence criterion itself, a biconditional like ours but with the operator side \#F = \#G dropped for the bare relation, and place their children at ``the initial stage of Russell--Frege's formal definition of cardinal integers''. Their children hold a weaker grade of the relation-without-the-operator that we will ascribe to the \Piraha{} (Section~\ref{sec:scopes}). Weaker in a precise sense: \citet{izardetal2014} decompose exact numerical equality into three principles, Identity, Addition/Subtraction, and Substitution; their 3-year-olds keep Identity, equality surviving a spatial rearrangement, and fail the other two on large sets.

\section{The puzzle}
Frege's Theorem (Section~\ref{sec:HP}) says HP entails the Peano axioms. The \Piraha{} grasp equinumerosity, in the enacted sense made precise below (Section~\ref{sec:about-Pirahas}). Yet they have no arithmetic. Where is the gap? We identify four cognitive loci: L0, the approximate number system (ANS); L1, the pairing relation (one-to-one matching); L2, the stable cardinal token ($= \#$); L3, the successor and generative counting. Possible confounding variables include general memory, motivation, task comprehension, and second-order quantification.

From the existing data:
\begin{enumerate}
\item \Piraha{} have an intact ANS (L0) and one-to-one matching (L1) in its enacted grade: exact within the glance untrained, out to ten where drilled (Section~\ref{sec:about-Pirahas}).
\item The failure is not plausibly reducible to motivation or comprehension: the failing tasks use the same materials, experimenters, and trial structure as the successful ones, thanks to the meticulous design of \citet{franketal2008}. \citet{everettandmadora2012} keep the same controls at the untrained village, where every participant passed the modeled practice trials at two and three. Nor is it general memory, in the sense of retaining absent objects: the orthogonal match fails with every spool in full view, while hidden-set responses show the ratio-scaled error characteristic of ANS, so an approximate summary \emph{was} stored and retrieved. The main boundary is therefore L1 to L2: enacted correspondence does not by itself yield a portable cardinal object. A neighboring dissociation, however, is also important: L2 is not the same as L3 (Table~\ref{tab:design}, phase III; Section~\ref{sec:buys}).\footnote{Consistent with this localization: exact symbolic number and approximate magnitude dissociate neurally \citep{dehaeneandcohen1997}; exact arithmetic is stored in a language-specific verbal format \citep{spelkeandtsivkin2001} yet survives the loss of grammar in agrammatic aphasia \citep{varleyetal2005}. The dependence is on the token, not on syntax, exactly where we place the boundary (hence, L2 but not L3).} What is never preserved is exact cardinal identity once direct pairing is blocked: the exact token is what is missing.
\item In Frank's battery the first-failing tasks require only \emph{storing and carrying an exact cardinality}, not successor, not iteration, not second-order quantification; at the untrained villages the failure arrives even earlier, at enacted pairing past the glance (Section~\ref{sec:about-Pirahas}). Either way the experiments fail before L3 is reached. And the gap is cardinality-specific: displacement in general is present in \Piraha{}. The hearsay evidential reports the unwitnessed, displacement by definition, and the verb system, at $\sim$65,000 forms, is anything but poor.
\item The equinumerosity relation is only an equivalence relation on concepts; it does not, on its own, deliver the cardinals as objects. That is the work of the term-former \#, and Frege's Theorem runs through \# (Section~\ref{sec:HP}), not through the relation alone. So even the full relation, lacking \#, would not yield the Peano axioms; the logical necessity of the operator matches the cognitive gap.
\end{enumerate}

The \Piraha{} have at most the right-hand side of HP (the relation); they never form the left-hand side (the operator). Throughout, `having the relation' means an enacted, finite correspondence ability: constructing a witnessing pairing on sets one can act on, exact for the few objects graspable at a glance and, untrained, no further. No grasp of the general second-order relation is claimed, and none is needed; by item 4, even the full relation would not be enough. The localization proceeds by elimination, \emph{neti neti} (`not this, not this'): not motivation, not comprehension, not general memory, not successor, not second-order quantification. The one locus the data cannot remove, and the minimal sufficient description of the deficit, is the missing term-former \#. L3 (successor, generative counting) is also absent, but the data never probe it; ``no successor'' and ``no recursion'' are explanatorily idle for these data. This localizes the gap by description, not by cause; Section~\ref{sec:scopes} addresses what the single-case design can and cannot support.

The four loci let us finish the metaphorical image of Section~\ref{sec:about-Pirahas}. Hands are L0 plus the glance-grade of L1: everyone's equipment, exact to about three, estimating above with that fixed CV. Threading is L1 in use: one thread from each balloon to its spool, a pairing witnessed on the table. Practice can stretch threading without adding anything else: at Xagiopai, Madora's coined words served as labeled thread, the threading was drilled, and it reached ten, always on the table. \emph{The stitch is \#: the move that closes a witnessed threading into one object that can leave the table.} A stitched bag is L2. A bag that grows as needed, a next size for every load, is L3, the generative count list. Hands give capacity; training can extend reach; the bag gives portability; the growing bag gives successor. 
In these terms, the data say something simple: hands are everywhere; threading can be stretched where it is locally drilled; no bag ever appears; and the question of growth never arises. In short: labels can stretch threading, but only a full count-list practice stitches a bag. One discipline is needed: the image tempts the stronger sentence, `they cannot stitch.' The data license only the description, `no bag was ever stitched.'

\section{Number words are the cognitive realization of Frege's \#}\label{sec:number-words}

Section~\ref{sec:HP} prompts the question: what does the operator \# exactly do? 
Following \citet{heck1999}, \# takes a concept and returns a stable \emph{object} (a named cardinal); in his phrase, ``number-words are not predicates of predicates, but proper names.'' This is a bootstrap of creating an object, not just an abstract idea. That is, \# produces something that can be named, stored, and re-identified. In Principle~\ref{HP}, the relation (RHS) lives in the witnessing arrangement of the sets; the operator (LHS) frees the cardinality from that arrangement.
From a cognitive perspective, \textbf{now, we can take the number word, store it in our mind, transport it, compare it across time, memory and modality.} 

From the study of \Piraha, \citet{franketal2008} reported that their participants pair off sets laid side by side (the relation) but cannot carry the result (no operator). \citeauthor{franketal2008} call what the \Piraha{} lack \emph{cognitive technology}, a suggestive but philosophically informal phrase.

So, the number words are the cognitive realization of \#, and that is what \emph{cognitive technology} was reaching for. Said with full precision: \# is the operator, and each number word names one of its values. What realizes the operator is the count-list practice as a whole. Stated as a definition: a practice realizes \# over a domain of concepts when it assigns to each a stable, re-identifiable token, and gives any two the same token if and only if they are equinumerous. The second clause is HP itself, read at the level of tokens; a practice that meets it does in behavior what \# does in the formalism. A bare label does not yet realize \#: it must be embedded in a routine that keeps it re-identifiable across absence, rearrangement, and delay, so that the token, not the scene, carries the cardinality. That is why Madora's labels extended enacted matching at Xagiopai without installing the operator (Section~\ref{sec:about-Pirahas}). ``Number words realize \#'' is our shorthand for that practice-level fact, and words are the canonical human carrier of the practice, not the only possible one (Section~\ref{sec:buys}). Frege's Theorem then reads: logical sufficiency (HP, in SOL, entails PA) with cognitive insufficiency (the mind does not get \# from perception alone), and number words close the gap. We now develop the identification along four lines: displacement, ostension, the second-order question, and psychologism.

\paragraph{Displacement:}
Apart from theoretical philosophy, this relates directly to \citet{hockett1960}'s displacement, with one refinement the data force. In Frank's five tasks, the two the \Piraha{} pass are exactly the two where the pairing can be \emph{enacted}, each balloon laid against its spool. The three they fail are the three where enactment is blocked and the cardinality itself must travel: across a rotation (orthogonal match, everything in full view), across a screen (hidden match), or across time (nuts-in-a-can). \citeauthor{franketal2008} describe the failing tasks in just these terms, exact quantity transferred across space and time. That is exactly where the L1/L2 boundary should appear. 
The relation needs only an arrangement one can act on; 
the operator is what frees the cardinality when that arrangement is hidden, moved, or reoriented.
So performance collapses precisely where a carried token is required, not where the sets leave sight. And our metaphorical bag of Section~\ref{sec:about-Pirahas} is no ordinary bag: an ordinary one carries the apples; a numeral carries the how-many and leaves the apples behind. Displacement in time is the dramatic case of the same demand, and the \Piraha{} do not lack displacement in general; they lack specifically a device that displaces exact cardinality. This gradient is Frank's; \citet{everettandmadora2012} place it at his one number-word-trained village (the local-training contrast of Section~\ref{sec:about-Pirahas}), and the same token-dependence of exact cardinality recurs in numerate adults under interference and in young children, so the cut does not rest on the \Piraha{} alone.

\paragraph{Ostension asymmetry}  Frank's participants absorbed the matching rule from a short demonstration and applied it essentially without error to 10 objects; at the untrained village the same modeling was understood at two and three and did not extend (Section~\ref{sec:about-Pirahas}). So showing conveys the relation where practice has prepared it. Nothing in the record conveys \#: not the demonstration, not Madora's months of drilling, not her words; the same participants who extend the relation to ten fail wherever the cardinal must be carried. The asymmetry is measured, not stipulated. And it is the cognitive echo of \emph{Grundlagen} \S 62 \citep{frege1884}: a cardinal cannot be given by pointing, only through the identity conditions that HP supplies; one can point at paired sets, never at \#F.
\paragraph{Second Order Logic?} It is an obvious question why we argue they lack \# and not SOL. On this, we conclude that SOL is never tested on them, and whether they have it or not, they lack \#. 
More specifically, the RHS of HP is itself second-order ($\exists$R, R is a bijection); but exhibiting one witness R is not commanding the quantifier. And a witness is what they exhibit: a pairing on small enacted sets, and out to ten where the token practice supports it (Section~\ref{sec:about-Pirahas}). The empirical dissociation isolates \#: they build the witness bijection (the relation in action) yet cannot name or carry the cardinal. They thread; they never stitch.

\paragraph{Psychologism and Mill} Frege wrote the \emph{Grundlagen} against two enemies at once: psychologism (never let the psychological contaminate the logical) and Mill's pebble-arithmetic empiricism \citep[\S 7]{frege1884}. A thesis like ours can sound like both, so we split what Mill collapsed. The paper describes the psychology of \emph{grasping} \#, not the nature of \# itself. \# stays a logical object throughout; \emph{we describe what a mind acquires} when it acquires the count lexicon. The operator serves as a model of the cognitive achievement, not a reduction of it; nothing about the truth or the a priori status of arithmetic is at stake here.

\section{What the identification buys us}\label{sec:buys}
We claim that this identification is generative for subsequent thoughts, not just decorative. Concretely speaking, 
\begin{itemize}
    \item \textbf{Symbolization, not recursion:} 
    The \Piraha{} reveal that the first relevant break occurs before recursion: at portable exact-cardinal symbolization (\#, L2). The L2/L3 separation is then supplied by the developmental and cross-cultural evidence (cardinal-principle knowers before successor knowers; the Munduruk\'u's partial tokens without successor; Table~\ref{tab:design}).
    Iteration (the count list) vs. hierarchical recursion is a later distinction \citep{fitchetal2005}; the \Piraha{} fail before either is reached. The identification locates the cliff precisely: you need \# before you need successor \citep{cheungetal2017}, and you need successor before you need recursion. Recursion here means the generative rules that make the numeral list unbounded, not syntactic embedding. Accounts that place the core of the number faculty at recursion \citep{watumulletal2014, hauserandwatumull2017, hiraiwa2017} begin one layer above where these data show the break.
    \item \textbf{Developmental stages:} L0 $\to$ L1 $\to$ L2 $\to$ L3 is not just logical ordering; in the populations tested so far, children's development aligns with this sequence. Cardinal-principle knowers precede successor knowers by years~\citep{sarneckaandcarey2008, cheungetal2017}. The Munduruk\'u have numerals to five, but only the words for one and two are used exactly; the rest name approximate ranges (`one hand,' the word for five, was offered for five to nine dots), and exact arithmetic fails beyond four or five, just past the glance \citep{picaetal2004}. \citeauthor{picaetal2004} draw the moral themselves: ``the availability of number names, in itself, may not suffice''; the Munduruk\'u lack a counting routine. The reach of \# tracks the tokens that function, not the words on the list. Read across populations instead of within one child, this same ordering is a cline rather than a wall: the \Piraha{} at one end (no exact tokens), the 3-year-olds at a still weaker grade of the same relation-without-operator \citep{izardetal2014}, numerate adults at the other. The realization of the operator is graded because its vehicle is; \# reaches exactly as far as the token practice that carries it, so no all-or-nothing biological threshold is needed to explain the gap.
    \item \textbf{Multiple realizability of \#:}  \# is realized by any practice that assigns a portable (i.e., words are not necessary), re-identifiable token to a cardinality: count words, tallies, fingers-as-tally, knots (the Inca quipu recorded numbers in knotted cords; \citealp{ascherandascher1981}), etc. It's not ``words specifically are magic'' but ``the operator needs some vehicle and these people have none.'' The data fit this token-level claim tightly: the \Piraha{} do not tally either \citep{gordon2004, everett2005}, and the homesigners never use their ten available fingers as one (Section~\ref{sec:scopes}); both cases lack \emph{all} stable cardinal tokens, not merely words. Two standing objections thereby become instances rather than counterexamples. 
Material-engagement accounts, on which tallies, notches, fingers, or other external tokens do the work \citep{overmann2017}, describe another realizer of \#. They do not compete with the thesis, since the thesis is not that spoken words are necessary, but that some stable token practice is. Likewise, the Indigenous
Australian children studied by \citet{butterworthetal2008}, reported to handle small exact numerosities
without count words, threaten only a spoken-word-necessity claim we never make.
If the success is confined to small/glance numerosities, it belongs with the hands; if it is carried exact cardinality beyond the glance, the thesis predicts some stable nonverbal token practice (spatial, body-based, or material) doing the work. Carried exact cardinality with no such practice anywhere would refute the identification, not become another realizer.
Frege never cared which signs carry the operator. This practice-first grain has an old relative, arithmetic as a \emph{technique} in \citet{wittgenstein1956}; we flag the kinship and stay out of rule-following.
    \item  \textbf{The cognitive technology and the compression contrast:} \# is practically the `cognitive technology' framed by \citet{franketal2008}. They conclude that,  `languages give their users a new route for the efﬁcient encoding of experience'. Borrowing ideas from signal processing literature, we say, \Piraha{} quantifiers \emph{h\'oi}/\emph{ho\'i}/\emph{ba\'agiso} and the whole ANS system are \textbf{lossy codes}, just as color terms, a standing controversy in the Whorfian literature~\citep{kempandregier2012, zaslavskyetal2018, regierandxu2017}. On the other hand, a numeral is \textbf{lossless} with respect to the cardinality invariant, a sufficient statistic under equinumerosity. Numeral systems vary in how far they reach and how conventionalized they are \citep{zariquieyetal2025}, but not in where ``seven'' sits: wherever a language has an exact ``seven,'' it denotes the same cardinality, because the denotation is pinned by the equivalence relation, not chosen by the community (a loose `about seven' is loose only against that fixed point; nothing analogous anchors `reddish'). The \Piraha{} lack the token-former that makes a lossless one possible.
\end{itemize}

\section{Arguments in papers may be further than they appear}
Many important papers make claims similar or complementary to ours, and we clarify where we stand: 
\begin{itemize}
    \item \citet{gallistelandgelman1992, gallistelandgelman2000}: The claim is that verbal numerals inherit their meaning from the preexisting nonverbal magnitude system. Words are calibrated pointers into the ANS magnitude line. This reading does not match the \Piraha{} and Munduruk\'u findings above, nor the binomial slips of counting aloud \citep{cordesetal2001, franketal2012}: the exact tokens behave like a second representational system running in parallel with the ANS, not like calibrated pointers into it.
    \item \citet{ripsetal2008}: The paper claims that natural number cannot be bootstrapped from perception or object-tracking alone. We agree and say that perception delivers the relation (matching) and locate the unbridgeable step at the term-former \#. 
    \item \citet{careyandbarner2019}: The nearest neighbor. Their bootstrapping account builds integer meanings from the count list, set-based quantification, and one-to-one correspondence, the ingredients of HP's right-hand side, and it is the leading story of how the count list acquires its meanings. It remains an acquisition story: it does not name the formal object that the acquisition installs, and it does not use the theorem (HP plus SOL entails PA) that makes the \Piraha{} gap puzzling in the first place. We supply those two pieces and stay silent on acquisition itself (Section~\ref{sec:scopes}).
    \item Piantadosi corpus: In multiple articles \citep{pittetal2022, oshaughnessyetal2022, piantadosi2023}, Piantadosi and colleagues model the acquisition trajectory in rich computational detail. Our contribution is orthogonal, identifying the formal object that the trajectory installs.
    \item \citet{snyderetal2018}: There is a question whether Frege's foundation should be cardinal (via \#) or ordinal (via successor/ordering). We do not settle that question. Nonetheless, the \Piraha show matching (cardinal) without successor (ordinal); children acquire cardinal meaning before successor meaning \citep{sarneckaandcarey2008} 
    which makes us lean towards the cardinal route. However, whether the full number system is purely cardinal is not claimed by us. Whether the ordinal count list is necessary to install cardinal meaning is an acquisition question we do not adjudicate (Section~\ref{sec:scopes}). Our claim is constitutive: once installed, number words realize \#.
\end{itemize}

\section{Scopes and limits (and non-causality)}\label{sec:scopes}
Here, we will argue that our localization by description is not a retreat from the causal question and is possibly the strongest well-posed claim in the neighborhood. 

\paragraph{The frame: single-case dissociation logic:} The \Piraha are a single case in the neuropsychological sense \citep{caramazza1986}: they form one \emph{confounded culture-cluster}, dozens of subjects but one unit. 
A single case validly establishes that matching and counting are separable (an existence claim), and establishes nothing about whether language causes the gap (a population/causal claim); the Whorfian and culture-first readings overreach exactly here.

\paragraph{The confounded variable:} The \Piraha{} simultaneously lack (i) cardinal labels (\#), (ii) the iterative counting habit, (iii) the cultural practice of generalizing beyond immediate experience, and possibly (iv) clause subordination. These four variables are fully confounded, so none can be shown to cause the gap; the paper therefore localizes it by \emph{description}, never by cause. Lacks (ii)-(iv) are precisely the culture/language causal stories the paper declines to adjudicate; from the currently available data about \Piraha, no strong Whorfian claim can be made.

\paragraph{The constitutive claim is not linguistic relativity:} One might grant the causal caution above and still hear the identification, ``number words realize \#,'' as Whorfianism through the back door. It is not. Strong linguistic relativity says the language a community speaks shapes how its members think in general. Our claim is far narrower and structural: one representational device (a stable, re-identifiable cardinal token) is what realizes one operation (naming a cardinality as an object), and nothing else about \Piraha{} cognition follows from it. \emph{The vehicle need not even be linguistic (Section~\ref{sec:number-words}); a tally or a knot would serve}. What carries the operator is a token practice, not a grammar, so the thesis is orthogonal to whether language colors perception, memory, or reasoning at large. It also points the other way from determinism on the question of limits. If the \Piraha{} were barred from exact number by the language they speak, that would be strong relativity; but they are not barred. Xagiopai took up the labels, and exact matching extended with them. This is constitution, not influence: Whorf has language act on a prior, independent cognition, and here there is no prior portable exact-cardinal representation to act on. Acquiring the practice is acquiring that device, not one faculty reshaping another. The lexicon is not a cage but an acquirable technology: universal, in that any mind can take it up; cultural, in that it does not grow without transmission. One could, with a full count list, install \# and then successor and then arithmetic; that is our thesis run forward, not a threat to it. The boundary we place at symbolization is a seam in how number is built (\# before successor before recursion), not a wall inside the person.

\paragraph{Possible de-confounders:} 
Several avenues could in principle untangle the four lacks; Table~\ref{tab:design} lays them out, and they differ sharply in strength.
\begin{itemize}
\item \textbf{Homesigners} are the best observational strand: culture held numerate, lexicon absent, the gap intact (detailed below).
\item \textbf{Within-population lexicon contrasts} are the heart of the de-confounding, since they hold the culture fixed and move only the count-list practice. The single-case caution bars causal inference from the \Piraha{} taken as one confounded unit; these contrasts are what escape it.

Among deaf Nicaraguan signers, \citet{flahertyandsenghas2011} find that those with an ordered set of number signs track large exact quantities while the non-counters do not, the culture-lexicon confound gone and the data real. The within-\Piraha{} echo is the Xagiopai-versus-untrained contrast of Section~\ref{sec:about-Pirahas}: quasi-interventional, since Madora inserted the four-to-ten words. The contrast is measured, not recollected: the trained and untrained rates of Section~\ref{sec:about-Pirahas}. A third data point strengthens it: two visiting Xagiopai adults, tested by \citeauthor{everettandmadora2012} under the same protocol, were error-free throughout. What is anecdotal is only the attribution to the words, Madora's recollection that matching improved after the words and not after other tactics; and villages with task exposure but no words, \citeauthor{gordon2004}'s included, stayed at baseline \citep{everettandmadora2012}. The residual confound is Xagiopai's decade of outside contact.
\item \textbf{Verbal interference} \citep{franketal2012, nedergaardetal2023}: the same dissociation in numerate adults with the count routine blocked. It pins the run-time dependence cleanly, but is silent on what installs the operator in the first place.
\item \textbf{The cross-linguistic gradient} \citep{picaetal2004, zariquieyetal2025} is dose-response-shaped, but between-culture confounded.
\item \textbf{Developmental precedence} \citep{sarneckaandcarey2008, cheungetal2017} is consistent with the ordering, but not decisive about cause.
\item \textbf{The ideal quasi-experiment} cannot be run: a sister language of the same Mura family with a count list in place would be it, but \Piraha{} is the family's last surviving member, and its one recorded relative, Bohur\'a, survives as roughly a hundred wordlist items \citep{everett1986}.
\end{itemize}
To claim something causal one needs an intervention, and the one intervention that would settle it, `a count lexicon without the counting practice,' is not well defined; so, strictly speaking, no well-defined causal effect of the words, as distinct from the practice, exists to estimate.

\begin{table}[htbp]
\centering
\small
\caption{Ten populations on the cognitive ladder, grouped by structural locus and ordered by ascending grade. $\checkmark$ = present; $\times$ = absent. L0 (the approximate number system) is intact in every row and omitted. Anum.\ = anumeric; num.\ = numerate. The count-routine and L-columns describe the individual, not the surrounding culture.}
\label{tab:design}
\resizebox{\textwidth}{!}{%
\begin{tabular}{lcccccr}
\toprule
Population & Culture & Own count routine & L1 (1-to-1) & L2 (\#) & L3 (succ.) & What it establishes \\
\midrule
\multicolumn{7}{l}{\textbf{I. Before portable cardinal objecthood}} \\
\quad 3-year-olds & num. & emerging & Identity & $\times$ & $\times$ & sameness $\neq$ enacted matching \\
\quad \Piraha{} untrained & anum. & $\times$ & $\checkmark$ ($\sim$3) & $\times$ & $\times$ & matching does not carry \\
\quad Homesigners & num. & $\times$ & $\checkmark$ ($\sim$3) & $\times$ & $\times$ & numerate culture $\neq$ token practice \\
\quad Nicaraguan non-counters & num. & $\times$ & $\checkmark$ ($\sim$3) & $\times$ & $\times$ & ambient counting $\neq$ own count routine \\
\quad \Piraha{} Xagiopai & anum. & local training; no count routine & $\checkmark$ (to 10) & $\times$ & $\times$ & trained L1 $\neq$ carried \# \\
\addlinespace
\multicolumn{7}{l}{\textbf{II. Partial or disrupted token support}} \\
\quad Numerate adults under interference & num. & access blocked & $\checkmark$ & blocked & blocked & blocking rehearsal blocks access to \# \\
\quad Munduruk\'u & partly num. & partial & $\checkmark$ ($\sim$5) & partial & $\times$ & partial tokens, partial \# \\
\addlinespace
\multicolumn{7}{l}{\textbf{III. Cardinal objecthood before successor}} \\
\quad CP-knowers ($\sim$4 to 5\,yr) & num. & $\checkmark$ & $\checkmark$ & $\checkmark$ & $\times$ & \# precedes successor \\
\addlinespace
\multicolumn{7}{l}{\textbf{IV. Integrated exact-number practice}} \\
\quad Nicaraguan counters & num. & $\checkmark$ & $\checkmark$ & $\checkmark$ & $\checkmark$ & count routine completes the ladder \\
\quad Numerate adults & num. & $\checkmark$ & $\checkmark$ & $\checkmark$ & $\checkmark$ & baseline: full ladder \\
\bottomrule
\end{tabular}%
}

\vspace{4pt}
{\footnotesize \textit{Notes:} ``Identity'' = cardinal equality survives rearrangement but not substitution or addition \citep{izardetal2014}; weaker than enacted matching. 
``Local training; no count routine'' = months of one-to-one matching practice plus coined labels for 4--10 at Xagiopai; the record supports extended enacted matching, but not a productive routine for carrying exact cardinality across hidden or orthogonal tasks.
CP-knowers = children who pass Give-N (cardinal meaning acquired) but have not yet grasped successor \citep{sarneckaandcarey2008, cheungetal2017}. ``Access blocked'' = verbal rehearsal suppressed; exact number drops toward approximate levels while access to the ordinary token routine is blocked \citep{franketal2012}.

\textit{Sources:} \Piraha{} untrained: \citet{gordon2004, everettandmadora2012}; Xagiopai: \citet{franketal2008, everettandmadora2012}; homesigners: \citet{spaepenetal2011}; Nicaraguan signers: \citet{flahertyandsenghas2011}; Munduruk\'u: \citet{picaetal2004}; CP-knowers: \citet{sarneckaandcarey2008, cheungetal2017}; numerate adults: \citet{franketal2012}.}
\end{table}

\paragraph{Homesigners} If we only work with \Piraha{} data, the culture and language are confounded. The Nicaraguan homesigners \citep{spaepenetal2011} separate this bundle. They are four deaf adults born into hearing families; they never learned any sign language and invented their own gestures at home; they live in a numerate culture, hold jobs, and handle money, yet their gesture system has no count list.

Across tasks, all four produced finger counts that tracked the target only approximately for larger sets, even when given unlimited time; the emblem from the introduction, the homesigner extending nine fingers for ten sheep with ten fingers available to tally, is \citeauthor{spaepenetal2011}'s own example. So the culture is numerate, the words are absent, and the gap persists: they hold the relation, in its enacted grade, and never form the quotient. 

This is a second, cleaner instance of HP-without-\#, and it does real subtractive work. Between the two cases the background inverts and the deficit stays. The homesigners generalize beyond the immediate every day, so neither lack (iii) above nor an anumeric surrounding can be what matters. The one lack the two cases share is the token practice itself. The subtraction stops there by necessity, not timidity. The labels and the routine are not two separable causes; they are components of one practice, and separating them is not a well-defined intervention, as noted above. And the nine-fingers emblem is the practice point in miniature: the tokens are literally in hand, ten of them; what is missing is the routine that would make them count. 

\paragraph{8-month teaching failure:} Eight months of Portuguese number classes, at the \Piraha's request, taught no one to count or add \citep{everett2005}. We do not explain it, and no existing account settles it. But the thesis survives, because it is constitutive, not interventional: a full counting routine cannot simply be handed over (thread can be; the stitch cannot). The asymmetry is structural: an axe maps onto motor schemas already in place, but each numeral has its identity only through its position in the routine, so the system cannot arrive one token at a time. The record even holds a second attempt. Keren Madora's Xagiopai sessions were arithmetic lessons too, again at the \Piraha's request \citep{everettandmadora2012}, and the arithmetic did not take there either. What took sits one rung lower: her coined words and drilling stretched exact \emph{matching}, while carrying stayed at baseline and the count list stayed absent (Section~\ref{sec:about-Pirahas}). The words arrived; \# did not. The only controlled evidence here is still \emph{removal} \citep{franketal2012}; no controlled installation has succeeded.\footnote{A cross-species echo of the installation failure: Matsuzawa's chimpanzee Ai (named for the Japanese word for love, \emph{ai}) acquired numeral-like labels only through years of explicit training, one numeral at a time, with none of the acceleration that marks the human cardinal-principle click \citep{matsuzawa1985, matsuzawa2009}. Tokens can be handed over; the routine that makes them realize \# does not come with them.} The candidate rescues are exactly the confounded causal stories this section declines to adjudicate: a childhood acquisition window, cultural resistance, or the possibility that the words can only arrive together with the whole counting routine. We name them and leave them.

\paragraph{What if?} Suppose one rejects \citet{franketal2008} outright and reads \citet{everettandmadora2012} as showing the \Piraha{} never match exactly beyond three. Section~\ref{sec:about-Pirahas} has already reconciled the two studies as one local-training contrast, and \citeauthor{everettandmadora2012} themselves conclude that exact recognition beyond three requires number terminology. That is token-dependence one rung below where our thesis places \#; support, not refutation. But even the harshest reading leaves the localization standing, because the clean L1-present, L2-absent dissociation never rested on the \Piraha{} alone: numerate adults under verbal interference \citep{franketal2012}, the homesigners \citep{spaepenetal2011}, and \citet{izardetal2014}'s 3-year-olds each build the pairing yet never form the token. 

\paragraph{What this does NOT prove:} No causal claim, no acquisition claim, no neural claim beyond the cited dissociation is stated here.

\section{Conclusion}
In this article, we conclude that equinumerosity is only the right-hand side of HP; arithmetic needs the left-hand side, and the left-hand side is the operator \#. Tracing back, we found that number words are the cognitive realization of Frege's \#. The \Piraha{} are essentially Hume's Principle without \#: the relation without the quotient. 
Furthermore, from Hume's principle to arithmetic the logical distance is zero: no further axiom is needed (Frege's Theorem). The cognitive distance that matters lies earlier, between ``as many as'' and \#, and it is the entire width of a portable token practice (the \Piraha, the homesigners, interference, the children). We have also been explicit about scope: which claims these data can support at all, and which they cannot. 

Hume's remark sat in the \emph{Treatise} as an observation about equality. Frege showed it contains all of arithmetic. The \Piraha{} show that between the remark and the arithmetic stands not a further inference but a name: one, then another, then another...

\section*{Acknowledgements}

Intellectual debts go to Terence Tao, whose comprehensible yet rigorous real-analysis books sparked my interest in axiomatic mathematics during my bachelor years, as independent reading alongside coursework. 
For this paper specifically, thanks to Mahasweta Acharjee for a conversation over Christmas vacation which led me to read analytical philosophy and Whorfian literature more rigorously; to Pallavi Chowdhuri for introducing me to Hockett and related works; to Dr. Subhrajyoti Roy for his reviews on the mathematical logic of the paper; to Professor Ganesh Babulal, my neuroscience advisor, for encouraging me to execute the idea into an article. 

\paragraph{AI disclosure:} Drafts of this paper were adversarially reviewed with the assistance of Claude (Anthropic); all ideas and inferences are strictly my own.

\bibliography{references}
\end{document}